\documentclass[twocolumn,showpacs,preprintnumbers,amsmath,amssymb,nobibnotes,prl]{revtex4}

\usepackage{graphicx}
\usepackage{dcolumn}
\usepackage{bm}

\begin{document}


\title{Analysis of Firing Behaviors in Networks of Pulse-Coupled Oscillators\\
with Delayed Excitatory Coupling}

\author{Wei Wu}
 \email{051018023@fudan.edu.cn}
\author{Tianping Chen}%
 \altaffiliation{Corresponding author}
 \email{tchen@fudan.edu.cn}
\affiliation{Key Laboratory of Nonlinear Mathematics Science, School
of Mathematical Sciences, Fudan University, Shanghai 200433,
P.~R.~China}



\date{\today}

\begin{abstract}
For networks of pulse-coupled oscillators with delayed excitatory
coupling, we analyze the firing behaviors depending on coupling
strength and transmission delay. The parameter space consisting of
strength and delay is partitioned into two regions. For one region,
we derive a low bound of interspike intervals, from which three
firing properties are obtained. However, this bound and these
properties would no longer hold for another region. Finally, we show
the different synchronization behaviors for networks with parameters
in the two regions.
\end{abstract}

\maketitle

For decades, complex networks have been focused on by scientists
from various fields, for instance, sociology, biology, chemistry and
physics, etc. In particular, networks of pulse-coupled oscillators,
as an important class of interconnected dynamical systems, have
gained increasing attentions because of their intimate relationship
to natural systems as diverse as cardiac pacemaker cells, flashing
fireflies, chirping crickets, biological neural networks, and
earthquakes (cf.
\cite{relationship-to-natual-systems,Buck1988,Mirollo1990}). A
pioneering work on modeling and analyzing pulse-coupled units was
done by Mirollo and Strogatz \cite{Mirollo1990}. Inspired by
Peskin's model for self-synchronization of the cardiac pacemaker,
they proposed a pulse-coupled oscillator model with undelayed
excitatory coupling to explain the synchronization of huge
congregations of South East Asian fireflies. With the framework of
the Mirollo-Strogatz model, many theoretical and numerical results
on pulse-coupled networks have been obtained
\cite{Vreeswijk1993-Chen1994-Corral1995,Mathar1996-Goel2002,Ernst1995,
Ernst1998,Timme2002-1,Timme2002-2,Kim2004,Wu2007}.

Pulse-coupling is difficult to handle mathematically because it
introduces discontinuous behavior into the otherwise continuous
model and so stymies most of the standard mathematical techniques
\cite{Strogatz1983}. Particularly for delayed pulse-coupling, the
mathematical analysis of collective dynamics of networks becomes a
challenging problem. Past research experience indicates that some
underlying facts and assumptions about firing behaviors play a
crucial role in mathematical analysis
\cite{Mirollo1990,Mathar1996-Goel2002,Ernst1995,Ernst1998,Timme2002-1,Timme2002-2,Wu2007}.
For example, in \cite{Mirollo1990,Mathar1996-Goel2002},
synchronization was proved by making use of the fact that the firing
order of oscillators is always preserved for complete and undelayed
pulse-coupling; in \cite{Ernst1998}, an assumption about firing
times made the analysis easier by reducing the number of case
distinctions; in \cite{Wu2007}, the proof of desynchronization was
essentially due to a low bound of interspike intervals.

In this Letter, networks of pulse-coupled oscillators with delayed
excitatory coupling are studied. We analyze the firing behaviors
depending on coupling strength and transmission delay. The parameter
space consisting of strength and delay is partitioned into two
regions. For one region, we give a low bound of interspike
intervals. By using the bound, three firing properties are derived,
which would be very helpful for discussing synchronization of
networks and stability of periodic solutions. Unfortunately, these
properties no longer hold for another region. Furthermore, the
different synchronization behaviors for networks with parameters in
the two regions are presented.

We consider a system of $N$ identical oscillators which are
pulse-coupled in a delayed excitatory manner. As in
\cite{Timme2002-1}, the coupling structure is specified by the sets
$\mathrm{Pre}(i)$ of presynaptic oscillators that send pulses to
oscillator $i$, or the sets $\mathrm{Post}(i)$ of postsynaptic
oscillators that receive pulses from oscillator $i$. A phase
variable $\phi_i(t)\in[0,1]$ is used to characterize the state of
the oscillator $i$ at time $t$. In the case of no interaction, the
dynamics of $\phi_i$ is given by
\begin{eqnarray}
\mathrm{d}\phi_i(t)/\mathrm{d}t=1,\label{law1}
\end{eqnarray}
namely, the cycle period of the free oscillator is 1. When $\phi_i$
reaches the threshold $\phi_i=1$, the oscillator $i$ fires and
$\phi_i$ jumps back instantly to zero, after which the cycle
repeats. That is,
\begin{eqnarray}
\phi_i(t)=1\Rightarrow\phi_i(t^{+})=0.\label{law2}
\end{eqnarray}
Because of the transmission delay, the oscillators interact by the
following form of pulse-coupling: if oscillator $i$ fires at time
$t$, it emits a spike instantly; after a delay time $\tau$, the
spike reaches all postsynaptic oscillators $j\in\mathrm{Post}(i)$
and induces a phase jump according to
\begin{eqnarray}
\phi_j(t+\tau)=
f^{-1}(\min[1,f(\phi_j((t+\tau)^{-}))+\varepsilon_{ij}]),\label{law3}
\end{eqnarray}
where $\varepsilon_{ij}>0$ is the coupling strength from oscillator
$i$ to oscillator $j$; and the function $f$ is twice continuously
differentiable, monotonously increasing, $f'>0$, concave down,
$f''<0$, and satisfies $f(0)=0$, $f(1)=1$. For a more detailed
introduction of the model, see
\cite{Mirollo1990,Ernst1995,Ernst1998,Timme2002-1,Timme2002-2,Wu2007}.
In this Letter, we further assume the following: (i) The coupled
system starts at time $t=0$ with a set of initial phases
$0<\phi_i(0)\leq1$; (ii) there is no self-interaction, i.e.,
$i\not\in\mathrm{Pre}(i)$ for any oscillator $i$; (iii) $0<\tau<1$;
and (iv) the coupling strengths are normalized such that for all
oscillator $i$,
$\sum_{j\in\mathrm{Pre}(i)}\varepsilon_{ji}=\varepsilon$ with
$0<\varepsilon<1$.

We partition the parameter space
$\mathcal{A}=\{(\tau,\varepsilon)|\,0<\tau<1,0<\varepsilon<1\}$ into
two regions
\begin{eqnarray*}
&&\mathcal{A}_1=\{(\tau,\varepsilon)\in\mathcal{A}|\,f(\tau)+\varepsilon<1\},\\
&&\mathcal{A}_2=\{(\tau,\varepsilon)\in\mathcal{A}|\,f(\tau)+\varepsilon\geq1\}.
\end{eqnarray*}
First of all, we use ``proof by contradiction" to prove that no
oscillator can fire twice in a time window of length $\tau$, if
parameters $(\tau,\varepsilon)\in\mathcal{A}_1$. Let $t_1$ and $t_2$
with $t_1<t_2$ be two successive firing times of oscillator $i$.
Suppose $t_2-t_1\leq\tau$. We claim that if $t_2>\tau$, there must
be some oscillator $i'\in\mathrm{Pre}(i)$ which fires more than once
in the time interval $(t_1-\tau,t_2-\tau]\cap[0,\infty)$. In fact,
this comes from the monotony and concavity assumption of the
function $f$. Since $f'>0$ and $f''<0$, we have that for any
$0<\delta<1$, if $0\leq\theta_1<\theta_2\leq f^{-1}(1-\delta)$, then
$f^{-1}(f(\theta_1)+\delta)-\theta_1<f^{-1}(f(\theta_2)+\delta)-\theta_2$,
namely the property (A7) in \cite{Ernst1998}. It implies that in the
same circle, the later the spike arrives, the larger the induced
phase jump is \cite{Ernst1995,Ernst1998}. Therefore, if all the
presynaptic oscillators $j\in\mathrm{Pre}(i)$ fire at most once in
$(t_1-\tau,t_2-\tau]\cap[0,\infty)$, then in the time interval
$(t_1,t_2]$ the sum of the phase jumps of oscillator $i$ is not more
than $f^{-1}(f(t_2-t_1)+\varepsilon)-(t_2-t_1)$, i.e., the sum
reaches its maximum if all spikes arrive at time $t_2$
simultaneously. It means $\phi_i(t_2)\leq
f^{-1}(f(t_2-t_1)+\varepsilon)\leq f^{-1}(f(\tau)+\varepsilon)<1$,
which contradicts that oscillator $i$ fires at $t_2$. Thus, there
exists some oscillator $i'\in\mathrm{Pre}(i)$ firing more than once
in $(t_1-\tau,t_2-\tau]\cap[0,\infty)$. Let
$t_3,t_4\in(t_1-\tau,t_2-\tau]\cap[0,\infty)$ with $t_3<t_4$ be two
successive firing times of oscillator $i'$. From $t_2-t_1\leq\tau$,
it follows that $t_4-t_3\leq\tau$. Similarly as above, if
$t_4>\tau$, then there must be some oscillator
$i''\in\mathrm{Pre}(i')$ which fires more than once in the time
interval $(t_3-\tau,t_4-\tau]\cap[0,\infty)$. Repeating the
derivation leads to a finite sequence of pairs of firing times:
\begin{eqnarray}
\{t_1,t_2\}\rightarrow\{t_3,t_4\}\rightarrow\cdots\rightarrow\{t_{2n-1},t_{2n}\}
\label{sequence}
\end{eqnarray}
which satisfies
$[t_{2k+1},t_{2k+2}]\subseteq(t_{2k-1}-\tau,t_{2k}-\tau]\cap[0,\infty)$
for $k=1,\ldots,n-1$; $t_{2n}\leq\tau$, $t_{2k}>\tau$ for
$k=1,\ldots,n-1$; and each term of (\ref{sequence}) is two
successive firing times of some oscillator. Particularly, $t_{2n-1}$
and $t_{2n}$ are two successive firing times of some oscillator
$i_0$. However, similarly as the argument of $\phi_i(t_2)<1$,
according to $t_{2n}\leq\tau$, $f(\tau)+\varepsilon<1$ and the
assumption that the coupled system starts at time $t=0$, we can get
$\phi_{i_0}(t_{2n})<1$. It contradicts that oscillator $i_0$ fires
at $t_{2n}$. This contradiction comes from our hypothesis
$t_2-t_1\leq\tau$. For a more detailed proof, see \cite{remark1}. As
a consequence, we get

{\it Theorem 1:} If parameters $(\tau,\varepsilon)\in\mathcal{A}_1$,
all interspike intervals of each oscillator in the coupled system
must be longer than the delay time $\tau$.

Here and throughout, ``interspike interval" is referred to as the
time between two successive firing activities of an oscillator.
However, as opposed to Theorem 1, at each
$(\tau,\varepsilon)\in\mathcal{A}_2$, the coupled system has
solutions in which some interspike intervals are not longer than
$\tau$. The simplest example is that the oscillators with initial
phases $\phi_1(0)=\cdots=\phi_N(0)$ fire synchronously with a period
$t=\tau$, if $(\tau,\varepsilon)\in\mathcal{A}_2$. In the rest of
the Letter, one will see that this can cause significantly different
dynamical behaviors of the system at
$(\tau,\varepsilon)\in\mathcal{A}_1$ and at
$(\tau,\varepsilon)\in\mathcal{A}_2$, especially the different
firing behaviors. Before discussing the difference of firing
behaviors, let us give some definitions and notations. Denote
$t^i_m$ the time at which oscillator $i$ fires its $m$-th time.
Clearly, the firing times $t^i_m$, $i=1\ldots,N$, $m\geq1$, are
determined by initial phases. For a given set of initial phases
$[\phi_1(0),\ldots,\phi_N(0)]$, the solution
$[\phi_1(t),\ldots,\phi_N(t)]$ is said to be a period-$d$ solution
if there exist a $\Delta t_0>0$, and positive integers $M$ and $d$
such that the firing times of arbitrary oscillator $i$ satisfies
$t^i_{m+d}-t^i_{m}=\Delta t_0$ for all $m\geq M$. For a given set of
initial phases $[\phi_1(0),\ldots,\phi_N(0)]$, the solution
$[\phi_1(t),\ldots,\phi_N(t)]$ is said to be a completely
synchronized solution if there exists a $T\geq0$ such that the phase
variables of arbitrary oscillators $i$ and $j$ satisfy
$\varphi_i(t)=\varphi_j(t)$ for all $t\geq T$. For the convenience
of later use, we let $\varepsilon_{ij}=0$ for
$j\not\in\mathrm{Post}(i)$. Then, the phase jump (\ref{law3}) also
holds for $j\not\in\mathrm{Post}(i)$.

By using Theorem 1, we conclude that if
$(\tau,\varepsilon)\in\mathcal{A}_1$, any solution of the coupled
system possesses the following properties:

{\it Property 1:} For oscillators $i$ and $j$ satisfying
$\varepsilon_{ij}=\varepsilon_{ji}$ and
$\varepsilon_{ki}=\varepsilon_{kj}$ for all $k\in
K_{ij}:=\{1,2,\ldots,N\}\setminus\{i,j\}$, if $t^i_{m_i}\leq
t^j_{m_j}$, then $t^i_{m_i+1}\leq t^j_{m_j+1}$, i.e., the firing
order of $i$ and $j$ is always preserved.

{\it Property 2:} For oscillators $i$ and $j$ satisfying
$\varepsilon_{ij}=\varepsilon_{ji}$ and
$\varepsilon_{ki}=\varepsilon_{kj}$ for all $k\in K_{ij}$, if
$t^i_{m_i}=t^j_{m_j}$, then $\phi_i(t)=\phi_j(t)$ for all $t\geq
t^i_{m_i}$.

{\it Property 3:} If $[\phi_1(t),\ldots,\phi_N(t)]$ is a completely
synchronized solution, then it is a period-one solution.

In fact, Theorem 1 implies that in the case of
$(\tau,\varepsilon)\in\mathcal{A}_1$, a spike of oscillator $i$ must
reach oscillators $j\in\mathrm{Post}(i)$ before $i$ emits the next
spike. Thus, for oscillators $i$ and $j$ satisfying
$\varepsilon_{ij}=\varepsilon_{ji}$ and
$\varepsilon_{ki}=\varepsilon_{kj}$ for all $k\in K_{ij}$, the
instantaneous synchronization $t^i_{m_i}=t^j_{m_j}$ can lead to
$\phi_i(t)=\phi_j(t)$ for all $t\geq t^i_{m_i}$ (Property 2). For
the same reason, if $[\phi_1(t),\ldots,\phi_N(t)]$ is a completely
synchronized solution, then
$t^i_{m+1}-t^i_{m}=1-[f^{-1}(f(\tau)+\varepsilon)-\tau]$ for all $i$
and $t^i_{m}\geq T$. That is to say, any completely synchronized
solution is a period-one solution with the final interspike interval
being $1-[f^{-1}(f(\tau)+\varepsilon)-\tau]$ (Property 3).

Due to space limitations, here we prove Property 1 for the case of
$N=2$. By Property 2, we only need to prove that if
$t^1_{m_1}<t^2_{m_2}$, then $t^1_{m_1+1}\leq t^2_{m_2+1}$. The proof
is divided into four cases:

\noindent {\bf Case 1}: $t^2_{m_2}\geq t^1_{m_1+1}$.

In this case, we have $t^1_{m_1+1}\leq t^2_{m_2}<t^2_{m_2+1}$.

\noindent {\bf Case 2}: $t^1_{m_1}+\tau\leq t^2_{m_2}<t^1_{m_1+1}$.

In this case, since $(t^2_{m_2},t^1_{m_1+1}]\subset
(t^1_{m_1}+\tau,t^1_{m_1+1}+\tau)$, oscillator $2$ cannot receive
any spikes from oscillator $1$ in the time interval
$(t^2_{m_2},t^1_{m_1+1}]$. This, combined with
$0=\phi_2((t^2_{m_2})^{+})<\phi_1(t^2_{m_2})<1$, leads to
$\phi_2(t)<\phi_1(t)$ for all $t\in(t^2_{m_2},t^1_{m_1+1}]$. It
implies $t^1_{m_1+1}<t^2_{m_2+1}$.

\noindent {\bf Case 3}: $t^2_{m_2}<t^1_{m_1}+\tau$ and
$\phi_1(t^1_{m_1}+\tau)>\phi_2(t^1_{m_1}+\tau)$.

Since $(t^1_{m_1}+\tau,t^1_{m_1+1}]\subset
(t^1_{m_1}+\tau,t^1_{m_1+1}+\tau)$ and
$\phi_2(t^1_{m_1}+\tau)<\phi_1(t^1_{m_1}+\tau)$, similarly as Case 2
we can get $\phi_2(t)<\phi_1(t)$ for all
$t\in(t^1_{m_1}+\tau,t^1_{m_1+1}]$. It implies
$t^1_{m_1+1}<t^2_{m_2+1}$.

\noindent {\bf Case 4}: $t^2_{m_2}<t^1_{m_1}+\tau$ and
$\phi_1(t^1_{m_1}+\tau)\leq\phi_2(t^1_{m_1}+\tau)$.

Since $(t^2_{m_2},t^1_{m_1}+\tau)\subset(t^1_{m_1},t^1_{m_1}+\tau)$,
by Theorem 1 oscillator $2$ cannot receive any spikes from
oscillator $1$ in $(t^2_{m_2},t^1_{m_1}+\tau)$. This, combined with
$0=\phi_2((t^2_{m_2})^{+})<\phi_1(t^2_{m_2})<1$, leads to
$\phi_2((t^1_{m_1}+\tau)^{-})<\phi_1((t^1_{m_1}+\tau)^{-})$. Let
$f_0=f(\phi_1((t^1_{m_1}+\tau)^{-}))-f(\phi_2((t^1_{m_1}+\tau)^{-}))$.
Because the spike emitted by oscillator $1$ at $t^1_{m_1}$ reaches
oscillator $2$ at $t^1_{m_1}+\tau$, we have
$f(\phi_2(t^1_{m_1}+\tau))-f(\phi_1(t^1_{m_1}+\tau))=\varepsilon_{12}-f_0$.
It can be claimed that
$f(\phi_2((t^2_{m_2}+\tau)^{-}))-f(\phi_1((t^2_{m_2}+\tau)^{-}))<\varepsilon_{12}-f_0$.
Indeed, this comes from the property (A5) in \cite{Ernst1998}:
$f(\theta_2)-f(\theta_1)>f(\theta_2+\delta)-f(\theta_1+\delta)$, if
$\theta_1<\theta_2$ and $0<\delta\leq1-\theta_2$. Denoting $\Delta
t=t^2_{m_2}-t^1_{m_1}$, from the property (A5) in \cite{Ernst1998}
we get
$f(\phi_2((t^2_{m_2}+\tau)^{-}))-f(\phi_1((t^2_{m_2}+\tau)^{-}))
=f(\phi_2(t^1_{m_1}+\tau)+\Delta t)-f(\phi_1(t^1_{m_1}+\tau)+\Delta
t)
<f(\phi_2(t^1_{m_1}+\tau))-f(\phi_1(t^1_{m_1}+\tau))=\varepsilon_{12}-f_0$.
Because the spike emitted by oscillator $2$ at $t^2_{m_2}$ reaches
oscillator $1$ at $t^2_{m_2}+\tau$, we have
$f(\phi_1(t^2_{m_2}+\tau))=\min[1,f(\phi_1((t^2_{m_2}+\tau)^{-}))+\varepsilon_{21}]$.
So, if $f(\phi_1(t^2_{m_2}+\tau))=1$, then from Theorem 1 it follows
that $f(\phi_2(t^2_{m_2}+\tau))<1=f(\phi_1(t^2_{m_2}+\tau))$; if
$f(\phi_1(t^2_{m_2}+\tau))<1$, then from the above claim it follows
that $f(\phi_2(t^2_{m_2}+\tau))-f(\phi_1(t^2_{m_2}+\tau))=
f(\phi_2((t^2_{m_2}+\tau)^{-}))-f(\phi_1((t^2_{m_2}+\tau)^{-}))-\varepsilon_{21}
<\varepsilon_{12}-f_0-\varepsilon_{21}=-f_0<0$. It implies
$t^1_{m_1+1}<t^2_{m_2+1}$.

In fact, we proved that for the case of $N=2$, if
$t^1_{m_1}<t^2_{m_2}$, then $t^1_{m_1+1}<t^2_{m_2+1}$. For the case
of $N>2$, the proof is similar, and also can be divided into the
above four cases. The distinction is that when $N>2$,
$t^1_{m_1+1}=t^2_{m_2+1}$ may happen in Cases 2-4. This derives from
the fact that two oscillators are likely to desynchronize, while the
other oscillators try to synchronize them \cite{Ernst1998}.

\begin{figure}[!t]
\includegraphics[width=2.75in]{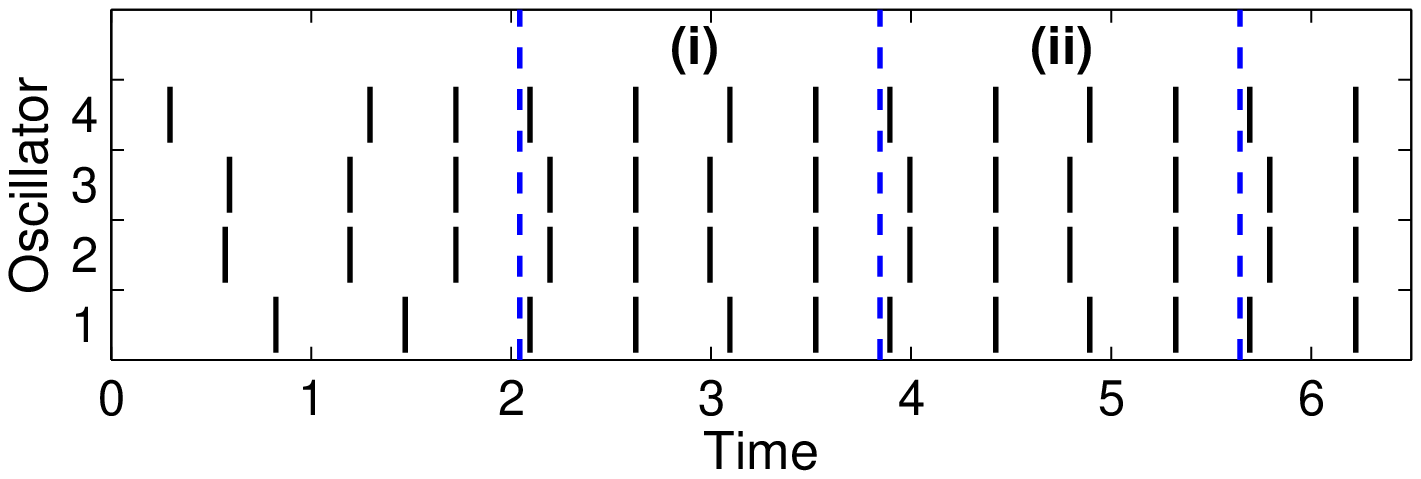}
\put(-3,0.8){\footnotesize (a)}

\vspace{1ex}
\includegraphics[width=2.75in]{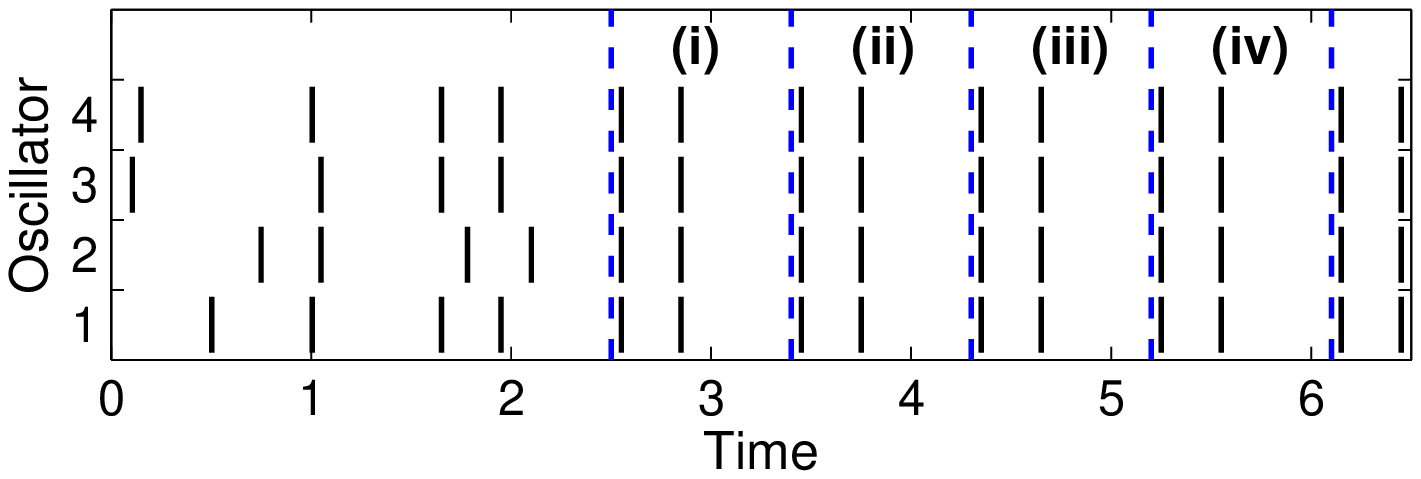}
\put(-3,0.8){\footnotesize (b)}

\caption{\label{fig:fig1} Firing times of four all-to-all
pulse-coupled oscillators with $\tau=0.9$ and $\varepsilon=0.6$.
Vertical dashed lines are used to indicate the boundaries of
periods. (a) Initial phases
$[\phi_1(0),\phi_2(0),\phi_3(0),\phi_4(0)]=[0.1766,0.4298,0.4079,0.7061]$.
Two periods (i) and (ii) are presented. (b) Initial phases
$[\phi_1(0),\phi_2(0),\phi_3(0),\phi_4(0)]=[0.4974,0.2492,0.8932,0.8501]$.
Four periods (i), (ii), (iii) and (iv) are presented.}
\end{figure}

Numerical analysis shows that from any initial phases, the coupled
system approaches a period solution with groups of synchronized
oscillators \cite{Ernst1995,Ernst1998,Timme2002-1,Wu2007,Kim2004}.
In larger networks, the oscillators can be divided into groups in a
combinatorial number of ways, and exponentially many periodic
solutions are present \cite{Timme2002-1}, which greatly increases
the complexity of firing behaviors. Properties 1-3 indicate that the
firing behaviors of the coupled system at
$(\tau,\varepsilon)\in\mathcal{A}_1$ are relatively simple. However,
when parameters $(\tau,\varepsilon)\in\mathcal{A}_2$, there may be
some solutions, which do not possess some or all of Properties 1-3.
It makes firing behaviors more complicated. Whether or not such
solutions exist depends on parameters $(\tau,\varepsilon)$ and
coupling strengths $\varepsilon_{ij}$. For the system with
$(\tau,\varepsilon)\in\mathcal{A}^0_2:
=\{(\tau,\varepsilon)\in\mathcal{A}|f(\tau)+\varepsilon>1\}$ and
$\varepsilon_{ij}=\varepsilon/(N-1)$, $i\neq j$ (hereinafter
referred to as all-to-all coupling), such solutions always exist.
Furthermore, for any such solution, there must be some interspike
intervals not exceeding the delay time $\tau$. Otherwise, by
previous arguments, the solution possesses Properties 1-3. As an
example, we simulate a network of $N=4$ all-to-all coupled
oscillators with $\tau=0.9$, $\varepsilon=0.6$. We use for $f$ an
example of the leaky integrate-and-fire model
\begin{eqnarray}
\frac{\mathrm{d}f(\phi)}{\mathrm{d}\phi}=-\bigg(\ln\frac{I}{I-1}\bigg)\cdot
f(\phi)+I\cdot\ln\frac{I}{I-1} \label{IF}
\end{eqnarray}
where $I=1.05$. In Fig. \ref{fig:fig1}(a), a period-four solution is
given. In this solution, the firing order of oscillators $1,2$ (or
$3,4$) is not preserved, e.g., $t^1_3<t^2_4$ but $t^1_5>t^2_6$; and
the instantaneous synchronization $t^i_{m_i}=t^j_{m_j}$ does not
mean $\phi_i(t)=\phi_j(t)$ for all $t\geq t^i_{m_i}$, e.g.,
$t^1_4=t^2_5$ but $t^1_5>t^2_6$. In Fig. \ref{fig:fig1}(b), a
period-two completely synchronized solution is given. In addition,
one can see that in Fig. \ref{fig:fig1}(a) and (b), most interspike
intervals of the oscillators are shorter than the delay time
$\tau=0.9$.

\setlength{\unitlength}{1in}

\setlength{\unitlength}{1in}
\begin{figure}
\includegraphics[width=2.25in]{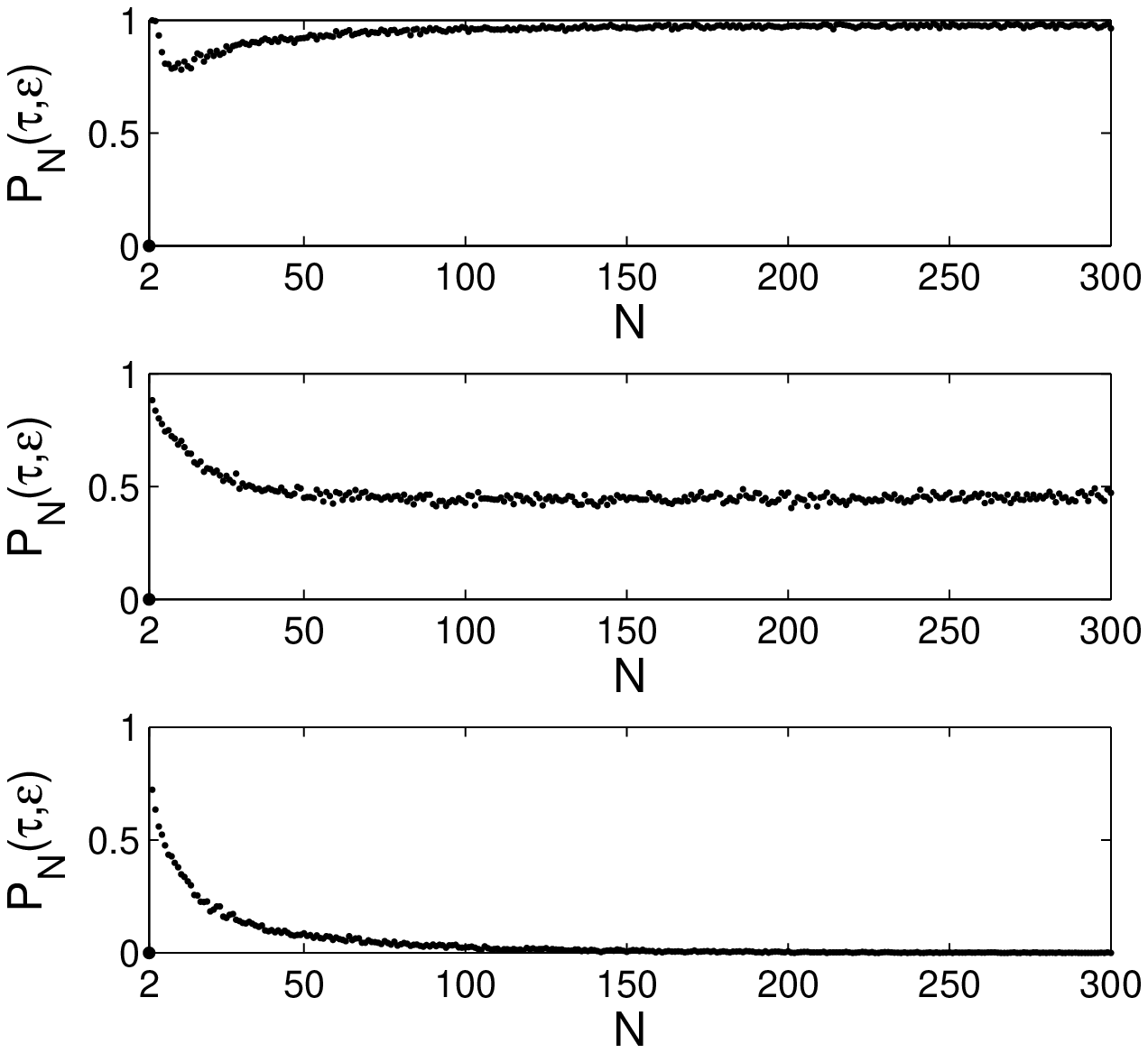}
\put(-2.5,1.85){\footnotesize(a)}
\put(-2.5,1.2){\footnotesize(b)}
\put(-2.5,0.55){\footnotesize(c)}

\vspace{-1ex}\caption{\label{fig:fig2} Dependence of
$P_N(\tau,\varepsilon)$ on network size $N$. (a) $\tau=0.55$,
$\varepsilon=0.4$. (b) $\tau=0.7$, $\varepsilon=0.35$. (c)
$\tau=0.8$, $\varepsilon=0.3$.}

\vspace{3ex}
\includegraphics[width=2.25in]{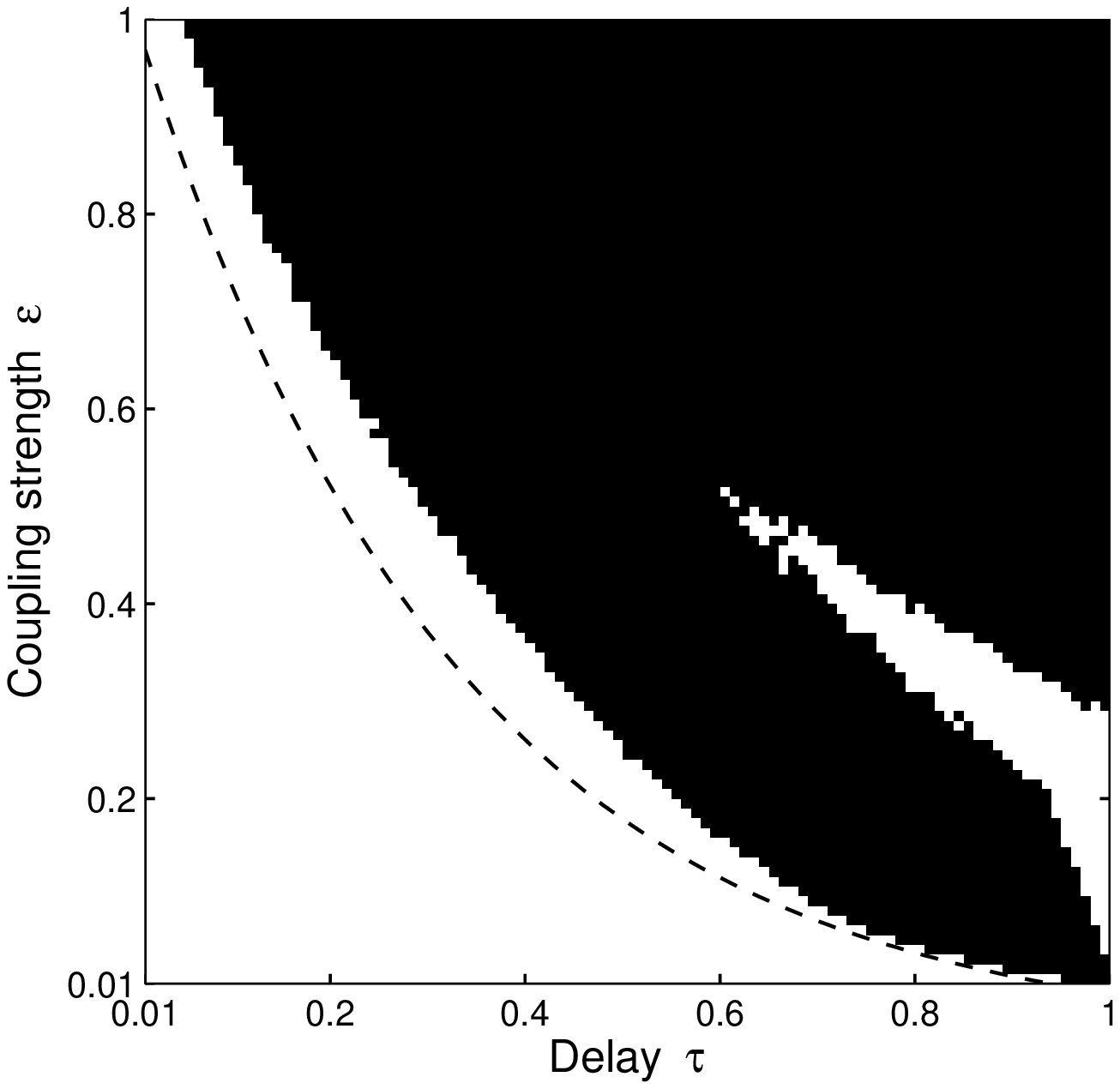}
\vspace{-1ex}\caption{\label{fig:fig3} Prevalence of completely
synchronized solutions for $(\tau,\varepsilon)\in\mathcal{A}_2$.
Parameters with $P_{100}(\tau,\varepsilon)>0$ are marked in black.
The dashed curve represents
$\{(\tau,\varepsilon)\in\mathcal{A}|f(\tau)+\varepsilon=1\}$.}
\end{figure}

Completely synchronized solutions, as a special type of periodic
solutions, have been widely studied
\cite{Mirollo1990,Mathar1996-Goel2002,Ernst1995,Ernst1998,Timme2002-2,Kim2004,Wu2007}.
The following analysis demonstrates the different synchronization
behaviors for networks with $(\tau,\varepsilon)\in\mathcal{A}_1$ and
$(\tau,\varepsilon)\in\mathcal{A}_2$. In \cite{Wu2007}, we proved
that under the assumption $f(2\tau)+\varepsilon<1$, from any initial
phases (other than $\phi_1(0)=\cdots=\phi_N(0)$), all-to-all
pulse-coupled oscillators with delayed excitatory coupling cannot
achieve complete synchronization. In fact, we can extend this result
to the case of $(\tau,\varepsilon)\in\mathcal{A}_1$ (see
\cite{remark2}). Interestingly, we found that when parameters
$(\tau,\varepsilon)\in\mathcal{A}_2$, completely synchronized
solutions become prevalent. In order to exhibit this, for networks
with all-to-all coupling, we numerically estimate the fraction
$P_N(\tau,\varepsilon)$ of the phase space
$\mathbf{\Phi}:=\{(\phi_1,\ldots,\phi_N)|0<\phi_i\leq1\}$ occupied
by initial phases of completely synchronized solutions. We still use
(\ref{IF}) for $f$. Fig. \ref{fig:fig2}(a)-(c) show the dependence
of $P_N(\tau,\varepsilon)$ on $N$ for $\tau=0.55$,
$\varepsilon=0.4$; $\tau=0.7$, $\varepsilon=0.35$; and $\tau=0.8$,
$\varepsilon=0.3$, respectively. More generally, we observed that
$P_N(\tau,\varepsilon)$ converges to a constant depending on
$(\tau,\varepsilon)$ as $N$ goes to infinity. For networks of
$N=100$, Fig. \ref{fig:fig3} shows the region of parameter space
$\mathcal{A}$ where completely synchronized solutions appear
($P_{100}(\tau,\varepsilon)>0$). For completely synchronized
solutions, there must be some interspike intervals not exceeding the
delay time $\tau$. Otherwise,  the system cannot be completely
synchronized (see \cite{Wu2007,remark2}). The performance of Fig.
\ref{fig:fig3} is supported by the observation \cite{Buck1988} of
flashing patterns of two firefly species {\it Photinus pyralis} and
{\it Pteroptyx malaccae}. For the species {\it P. pyralis}, the
normalized delay (neural delay/endogenous flashing period) is
$\approx0.03$. The whole group of the species rarely synchronizes
flashing; instead, wave, chain or sweeping synchrony has been
reported. For the species {\it P. malaccae}, the normalized delay is
$\approx0.36$, and perfect synchrony is usually achieved.

In summary, our analysis demonstrates  different dynamics for
pulse-coupled networks with $(\tau,\varepsilon)\in\mathcal{A}_1$ and
$(\tau,\varepsilon)\in\mathcal{A}_1$. For the region
$\mathcal{A}_1$, we derive a low bound of interspike intervals and
three firing properties, which provide a basis for future researches
addressing the dynamics in networks, e.g., stability of periodic
solutions. The difference of synchronization presented at the end of
the Letter is useful for understanding and interpreting
synchronization phenomena in some natural systems.

This work was supported by the National Science Foundation of China
under Grants 60574044 and 60774074.

\end{document}